# Generalized-Hypergeometric Solutions of the General Fuchsian Linear ODE Having Five Regular Singularities

**Artur Ishkhanyan** [1,2] **and Clemente Cesarano** [3,*]

1   Department of General Physics, Russian-Armenian University, Yerevan 0051, Armenia
2   Matter Wave Physics Department, Institute for Physical Research, Ashtarak0203, Armenia
3   Section of Mathematics-International Telematic University Uninettuno, C.so Vittorio Emanuele II, 39, 00186 Roma, Italy
*   Correspondence: c.cesarano@uninettunouniversity.net

**Abstract:** We show that a Fuchsian differential equation having five regular singular points admits solutions in terms of a single generalized hypergeometric function for infinitely many particular choices of equation parameters. Each solution assumes four restrictions imposed on the parameters: two of the singularities should have non-zero integer characteristic exponents and the accessory parameters should obey polynomial equations.

---

## 1. Introduction

During the last decades, the five Heun equations and their generalizations became a subject of intensive investigations in the context of numerous advanced problems of contemporary fundamental and applied research. Being natural generalizations of the equations of the hypergeometric class, these equations frequently appear, for instance, in contemporary classical and non-classical physics. Some representative examples include nuclear reactions, neutrino oscillations, quantum mechanics and quantum information theory, nanophysics, ultralow-temperature physics (Bose–Einstein condensation), fluid mechanics, statistical physics, astrophysics, cosmology, etc. (see, e.g., [1,2] and references therein)

The general Heun equation, from which the other equations originate, is the most general Fuchsian second-order linear differential equation having four regular singular points [3,4]. Here we consider the linear ordinary differential equation:

$$\frac{d^2u}{dz^2} + \left( \frac{\gamma}{z} + \frac{\delta}{z-1} + \frac{\varepsilon}{z-a} + \frac{\varepsilon_1}{z-a_1} \right) \frac{du}{dz} + \frac{\alpha\beta z^2 - \theta_1 z - \theta_0}{z(z-1)(z-a)(z-a_1)} u = 0 \tag{1}$$

which presents the most general Fuchsian equation having *five* regular singular points [3,4]. The singularities are located at $z = 0, 1, a, a_1$ and at $z = \infty$. The regularity of the singularity at infinity is ensured by the Fuchsian condition

$$1 + \alpha + \beta = \gamma + \delta + \varepsilon + \varepsilon_1 \tag{2}$$

If $\varepsilon_1 = 0$ and $\theta_1 = a_1\alpha\beta - \theta_0 / a_1$, Equation (1) becomes the general Heun equation [5]

$$\frac{d^2u}{dz^2} + \left( \frac{\gamma}{z} + \frac{\delta}{z-1} + \frac{\varepsilon}{z-a} \right) \frac{du}{dz} + \frac{\alpha\beta z - q}{z(z-1)(z-a)} u = 0 \tag{3}$$

with $q = -\theta_0 / a_1$, which is the most general Fuchsian equation having four regular singular points [6–8]. The parameter $q$ which is not related to the characteristic exponents of the singularities, is referred to as the accessory parameter. For Equation (1), there are two accessory parameters $-\theta_1$ and $\theta_0$.

There are many other choices of parameters that reduce Equation (1) to the general Heun Equation (3). For instance, an obvious choice is $\varepsilon = a = \theta_0 = 0$. This and the previous reduction cases can be made symmetric if the parameters $\theta_1$ and $\theta_0$ are expressed through new parameters $q$ and $q_1$ such that $\theta_1 = q + q_1$ and $\theta_0 = qq_1$. Then, for the triads $(a, \varepsilon, q) = (0, 0, 0)$ and $(a_1, \varepsilon_1, q_1) = (0, 0, 0)$ Equation (1) becomes the general Heun equation with the third finite singularity being located at $z = a_1$ and $z = a$, respectively. For this reason, in this case, one may conventionally think of the parameters $q$ and $q_1$ as accessory parameters associated with the singularities $a$ and $a_1$, respectively.

To find solutions of Equation (1) in terms of the generalized hypergeometric functions [9,10], we apply the approach developed in [11,12] for Heun-type equations which originate, via coalescence of singularities, from the general Heun equation [6–8]. A main result we report in the present paper is that there exist infinitely many particular choices of parameters for which Equation (1) having five irreducible regular singularities admits solutions in terms of a single generalized-hypergeometric function. We should keep in the mind the possible reductions of Equation (1) to the Heun equation in order to, in due course, identify the cases when the constructed solutions in fact apply to the general Heun equation, not to Equation (1) with irreducible singularities. These reducible cases are not taken into account.

## 2. The Approach

To construct the solutions, we examine the cases, when a power-series Frobenius solution [13] of Equation (1) reduces to the generalized hypergeometric series [9,10]

$$_rF_s\left(a_1, \ldots, a_r; b_1, \ldots, b_s; z\right) = \sum_{n=0}^{\infty} c_n z^n \tag{4}$$

for which the coefficients obey the two-term recurrence relation

$$\frac{c_n}{c_{n-1}} = \frac{1}{n} \frac{\prod_{k=1}^{r}\left(a_k - 1 + n\right)}{\prod_{k=1}^{s}\left(b_k - 1 + n\right)} \tag{5}$$

The Frobenius solution of Equation (1) for the vicinity of the regular singularity located at $z = 0$ is given as

$$u(z) = z^\mu \sum_{n=0}^{\infty} c_n z^n \tag{6}$$

where the characteristic exponent $\mu$ may adopt two values: $\mu = 0$ or $1 - \gamma$, and the expansion coefficients obey a *four*-term recurrence relation [3,4]:

$$R_n c_n + Q_{n-1} c_{n-1} + P_{n-2} c_{n-2} + S_{n-3} c_{n-3} = 0 \tag{7}$$

For the exponent $\mu = 0$ the coefficients are given as

$$R_n = aa_1 n(\gamma + n - 1) \tag{8}$$

$$Q_n = \theta_0 - \big((\gamma + n - 1)(a + a_1 + aa_1) + a\varepsilon_1 + a_1\varepsilon + aa_1\delta\big)n \tag{9}$$

$$P_n = \theta_1 + \big((n + \alpha + \beta)(1 + a + a_1) - \delta - a\varepsilon - a_1\varepsilon_1\big)n \tag{10}$$

$$S_n = -(n + \alpha)(n + \beta) \tag{11}$$

We note that the series (6) with such coefficients may terminate only if $\alpha$ or $\beta$ is a non-positive integer (besides, $\theta_0$ and $\theta_1$ should obey certain polynomial equations). Having this in the mind and based on the experience gained in treating the general and single-confluent Heun equations [11,12], we try the ansatz (5) for $r = N + 2$, $s = N + 1$, and

$$a_1,....,a_N, a_{N+1}, a_{N+2} = \alpha, \beta, 1 + e_1,...,1 + e_N \tag{12}$$

$$b_1,....,b_N, b_{N+1} = \gamma, e_1,...,e_N \tag{13}$$

This means that we look for a solution of Equation (1) as

$$u = {}_{N+2}F_{1+N}(\alpha, \beta, 1 + e_1,...,1 + e_N; \gamma, e_1,...,e_N; \sigma z) \tag{14}$$

where $\sigma$ is a scaling factor to be defined later. The series representation for this function is governed by a two term recurrence relation for coefficients explicitly written as

$$\frac{c_n}{c_{n-1}} = \sigma \frac{(\alpha - 1 + n)(\beta - 1 + n)}{(\gamma - 1 + n)n} \prod_{k=1}^{N} \frac{e_k + n}{e_k - 1 + n} = s_0 \frac{S_{n-1}}{R_n} \prod_{k=1}^{N} \frac{e_k + n}{e_k - 1 + n} \tag{15}$$

where $s_0 = -\sigma aa_1$. With this, the recurrence relation (7) is rewritten as

$$s_0^3 S_{n-1} S_{n-2} \prod_{k=1}^{N}(e_k + n) + s_0^2 Q_{n-1} S_{n-2} \prod_{k=1}^{N}(e_k - 1 + n) + \\ s_0 R_{n-1} P_{n-2} \prod_{k=1}^{N}(e_k - 2 + n) + R_{n-1} R_{n-2} \prod_{k=1}^{N}(e_k - 3 + n) = 0. \tag{16}$$

This is a polynomial equation in summation index $n$ of the degree utmost $n^{N+4}$. To ensure the fulfillment of this equation for all $n$, we demand the coefficients to identically vanish. For the coefficient of the highest order term we have $(a + s_0)(a_1 + s_0)(aa_1 + s_0) = 0$, hence,

$$s_0 = -a, -a_1, -aa_1 \tag{17}$$

The equation for the term proportional to $n^{N+3}$ reads $a^2(a - 1)a_1^2(a_1 - 1)(N + \varepsilon + \varepsilon_1) = 0$, hence, since $a, a_1 \neq 0, 1$,

$$\varepsilon + \varepsilon_1 = -N \tag{18}$$

Calculating the sum of all coefficients, we derive another relation:

$$\theta_0 = s_0 \alpha\beta \frac{\prod_{k=1}^{N}(1 + e_k)}{\prod_{k=1}^{N} e_k} \tag{19}$$

Finally, one more relation is derived if one checks the value adopted by the product $\prod_{k=0}^{N}(k + \varepsilon_1)$. Notably, the product turns to be zero. In the force of Equation (18), this means that $\varepsilon$ and $\varepsilon_1$ are *non-positive integers* such that $\varepsilon + \varepsilon_1 = -N$. With these observations, Equation (16) is reduced to a polynomial equation for $n$ of the degree $N$: $\sum_{m=0}^{N} A_m n^m = 0$, where $A_m$ are functions of the parameters of the starting Equation (1) and the parameters $\theta_1$ and $e_1,...,e_N$. Resolving then the equations $A_m = 0$, $m = 0,1,..,N$, we arrive at a polynomial equation of the degree $N+1$ for $\theta_1$ and the corresponding parameters $e_1,...,e_N$ for the generalized-hypergeometric solution (14) of Equation (1).

### 3. Explicit Solutions of Lowest Orders $N = 0,1,2$

Consider the choice $s_0 = -aa_1$, that is $\sigma = 1$.

The case $N = 0$ is almost trivial:

$$N = 0: \ \varepsilon = \varepsilon_1 = 0, \ \theta_1 = (a + a_1)\alpha\beta, \ \theta_0 = -aa_1\alpha\beta, \ u = {}_2F_1(\alpha,\beta;\gamma;z) \tag{20}$$

In this case the parameters $\theta_1$ and $\theta_0$ are such that the singularities at $z = a$ and $z = a_1$ both disappear so that Equation (1) is reduced to the ordinary Gauss hypergeometric equation.

For $N = 1$ we have two choices: $(\varepsilon,\varepsilon_1) = (-1,0)$ and $(\varepsilon,\varepsilon_1) = (0,-1)$, both reproducing the result for the general Heun equation [14–18,12]. Here is the explicit solution for the first choice:

$$N = 1: \ \varepsilon = -1, \ \varepsilon_1 = 0: \ u = {}_3F_2(\alpha,\beta,1+e_1;\gamma,e_1;z) \tag{21}$$

$$\theta_1 = a_1\alpha\beta - \theta_0 / a_1, \ \theta_0 = -aa_1\alpha\beta\frac{1+e_1}{e_1} \tag{22}$$

$$a(\alpha - e_1)(\beta - e_1) + (\gamma - 1 - e_1)e_1 = 0 \tag{23}$$

Because of the particular value of $\theta_1$, in this case the singularity at $z = a_1$ disappears so that Equation (1) is reduced to the general Heun equation. Resolving the second Equation (22) with respect to $e_1$ and substituting the result into (23), we arrive at a quadratic equation for $\theta_0$, which provides two values of this accessory parameter for which the solution (21) applies. This equation reads

$$\frac{\theta_0^2}{a_1^2} + \left(1 + a(\alpha + \beta + 2\alpha\beta) - \gamma\right)\frac{\theta_0}{a_1} + a\alpha\beta\left(a(1+\alpha)(1+\beta) - \gamma\right) = 0 \tag{24}$$

For $N = 2$ there are three possibilities: $(\varepsilon,\varepsilon_1) = (-2,0)$, $(\varepsilon,\varepsilon_1) = (0,-2)$, and $(\varepsilon,\varepsilon_1) = (-1,-1)$. For the first two cases a singularity disappears and the result reproduces that for the general Heun equation. For the first choice the solution reads [12,18]

$$N = 2: \ \varepsilon = -2, \ \varepsilon_1 = 0: \ u = {}_4F_3(\alpha,\beta,1+e_1,1+e_2;\gamma,e_1,e_2;z) \tag{25}$$

with parameters $e_1,e_2$ being defined by equations

$$\theta_0 = -aa_1\alpha\beta\frac{e_1+1}{e_1}\frac{e_2+1}{e_2}, \ a^2 = \frac{\prod_{k=1}^{k=2}(1-\gamma+e_k)e_k}{\prod_{k=1}^{k=2}(\alpha-e_k)(\beta-e_k)} \tag{26}$$

This solution applies if $\theta_1 = a_1\alpha\beta - \theta_0 / a_1$ and $\theta_0$ obeys the cubic equation

$$\begin{aligned}&\left(\omega^2 - \omega\left(\gamma - 2 + (1 - \alpha - \beta)a\right) + 2a(a-1)\alpha\beta\right)\times\\&\left(\omega + 2a(1+\alpha+\beta) + 2 - 2\gamma\right) + 2\omega a(a-1)(1+\alpha)(1+\beta) = 0,\end{aligned} \tag{27}$$

where $\omega = a\alpha\beta + \theta_0 / a_1$.

The third choice $(\varepsilon,\varepsilon_1) = (-1,-1)$ suggests a non-trivial solution again in terms of the functions ${}_4F_3$:

$$N = 2: \ \varepsilon = -1, \ \varepsilon_1 = -1: \ u = {}_4F_3(\alpha,\beta,1+e_1,1+e_2;\gamma,e_1,e_2;z) \tag{28}$$

$$\theta_1 = \alpha\left((1 + a + a_1)\beta - (\delta - a - a_1) + \frac{(1 + \alpha - \gamma)\prod_{k=1}^{k=2}(1 + \alpha - e_k)}{\prod_{k=1}^{k=2}(\alpha - e_k)}\right) \tag{29}$$

$$\theta_0 = -aa_1\alpha\beta\frac{1 + e_1}{e_1}\frac{1 + e_2}{e_2} \tag{30}$$

$$(\alpha - \beta)(\delta - a - a_1) = \frac{\alpha(1 + \alpha - \gamma)\prod_{k=1}^{k=2}(1 + \alpha - e_k)}{\prod_{k=1}^{k=2}(\alpha - e_k)} - \frac{\beta(1 + \beta - \gamma)\prod_{k=1}^{k=2}(1 + \beta - e_k)}{\prod_{k=1}^{k=2}(\beta - e_k)} \tag{31}$$

$$aa_1 = \frac{\prod_{k=1}^{k=2}(1 - \gamma + e_k)e_k}{\prod_{k=1}^{k=2}(\alpha - e_k)(\beta - e_k)} \tag{32}$$

where Equations (29) and (30) parameterize $\theta_0, \theta_1$ through the parameters $e_1$ and $e_2$, while the last two equations define the latter parameters. Eliminating $e_1, e_2$ from Equations (29) and (30), it can be shown that $\theta_1$ can be presented as a rational function of $\theta_0$:

$$\begin{aligned}
\theta_1 = &\ 2 - \gamma + (1 + \alpha)\beta(a + a_1) - \\
&\left[aa_1\left[\alpha\left(2a(1 + \alpha)(1 + \beta)(2 + \beta)a_1 + \gamma^2(a + a_1) - (2 + \beta)\gamma\left(a_1 + a\left(1 + (1 + \alpha + \beta)a_1\right)\right)\right) + \right.\right. \\
&\left.\left.\left(2(1 + \alpha)(1 + \beta) - (1 + \alpha + \beta)\gamma\right)\theta_0\right]\right] / \\
&\left[aa_1\alpha(1 + \beta)\left(a(1 + \alpha)(2 + \beta)a_1 - \gamma(a + a_1)\right) + \right. \\
&\left.\left(a_1 - \gamma(a + a_1) + a\left(1 + (1 + \beta + \alpha(3 + 2\beta))a_1\right)\right)\theta_0 + \theta_0^2\right],
\end{aligned} \tag{33}$$

while $\theta_0$ obeys a (rather cumbersome) *fourth*-order polynomial equation the coefficients of which are functions of the parameters $a, a_1, \alpha, \gamma, \delta$:

$$\theta_0^4 + A_3\theta_0^3 + A_2\theta_0^2 + A_1\theta_0 + A_0 = 0 \tag{34}$$

with

$$A_3 = 2\left(aa_1(2\alpha\beta + 2\alpha + 2\beta + 1) - (a + a_1)(\gamma - 1)\right) \tag{35}$$

$$\begin{aligned}
A_2 = &\ a_1^2 a^2\left((6\alpha(\alpha + 2) + 5)\beta^2 + 4\beta(3\alpha + 5)\beta + \alpha(5\alpha + 6) + 6\beta + 1\right) - \\
&\ a_1 a^2\gamma\left(\alpha(6\beta + 5) + 5\beta + 2\right) + a_1 a^2\left(4\alpha(\beta + 1) + 4\beta + 2\right) + \left(a^2 + a_1^2\right)(\gamma - 1)^2 + \\
&\ a_1 a\left(-a_1\gamma(6\alpha\beta + 5\alpha + 5\beta + 2) + 2a_1(2\alpha\beta + 2\alpha + 2\beta + 1) + 3(\gamma - 2)\gamma + 2\right),
\end{aligned} \tag{36}$$

$$\begin{aligned}
A_1 = &\ aa_1\left(a_1(a_1(2(\alpha + 1)a^2(\beta + 1))\left((2\alpha(\alpha + 2) + 1)\beta^2 + \alpha(4\alpha + 7)\beta + \alpha(\alpha + 1) + \beta\right) - \right. \\
&\ a\gamma\left((2\alpha(3\alpha + 5) + 3)\beta^2 + (2\alpha(5\alpha + 7) + 3)\beta + 3\alpha(\alpha + 1)\right) + \\
&\ 2(\alpha + 1)a(\beta + 1)(\alpha\beta + \alpha + \beta) + (\gamma - 1)\gamma(2\alpha\beta + \alpha + \beta)) - \\
&\ \gamma\left(a\left(a\alpha^2(2\beta(3\beta + 5) + 3) + \alpha\left(a(2\beta(5\beta + 7) + 3) + 8(\beta + 1)\right) + \beta(3a(\beta + 1) + 8) + 6\right) + 2 \\
&\ \gamma^2\left(a(6\alpha\beta + 4\alpha + 4\beta + 2) + 3\right) + 2a(\alpha + 1)(\beta + 1)\left(a(\alpha\beta + \alpha + \beta) + 2\right) - \gamma^3\right) + \\
&\ a(\gamma - 1)\gamma\left(a(2\alpha\beta + \alpha + \beta) - \gamma + 2\right)),
\end{aligned} \tag{37}$$

$$
\begin{aligned}
A_0 = {} & a^2 a_1^2 \alpha\beta(a\gamma^2\left(a\left(\alpha+1\right)\left(\beta+1\right)-\gamma+1\right)+ \\
& a_1\big(\gamma\big(a\gamma\big(3\alpha\beta+4\alpha+4\beta+4\big)-a\left(\alpha+1\right)\left(\beta+1\right)\big(a\left(2\alpha\beta+3\alpha+3\beta+4\right)+2\big)-\gamma^2+\gamma\big)+ \\
& -a_1\left(\alpha+1\right)\left(\beta+1\right)\big(a\left(\alpha+2\right)\left(\beta+1\right)-\gamma\big)\big(a\left(\alpha+1\right)\left(\beta+2\right)-\gamma\big)\big)\big).
\end{aligned}
\tag{38}
$$

This result has been checked by Mathematica.

## 4. Discussion

Thus, we have shown that there exist infinitely many particular choices of parameters for which a Fuchsian differential equation having five regular singular points admits solutions in terms of a single generalized-hypergeometric function. This is the main result of this paper.

It should be noted that there have been several studies towards the solution of Fuchsian equations with more than four singularities (these are more complicated equations than the general Heun equation) in terms of simpler mathematical functions of the hypergeometric class (see, e.g., [18–24]). In particular, based on Scheffé's lemma [19], several such solutions have been derived in [20,21]. A computer algebra solver to find solutions having rational function arguments has been presented in [22] (note that the solver works with equations for which the coefficients do not involve continuously variable parameters).

The results we have reported here, that is the closed-form explicit solutions for equations having five regular singular points, are in line with similar results derived for the Heun-type equations [11,12]. The results are in agreement with the conjecture by Takemura [18] which suggests that generalized-hypergeometric solutions exist for any Fuchsian differential equation having $3 + M$ regular singularities of which $M$ are apparent [25]. As we have seen, in our case each generalized-hypergeometric solution assumes four restrictions imposed on the equation parameters, namely, two of the singularities should have non-zero integer characteristic exponents and two accessory parameters should obey certain polynomial equations. It can be shown that the latter restrictions warrant that the two singularities with non-zero integer characteristic exponents are indeed apparent.

Fuchsian linear ordinary differential equations having five regular singularities appear in many branches of contemporary physics and mathematics research (see, e.g., [1,26,27]). Based on the experience gained in application of such solutions in the case of the Heun-type equations (see, e.g., [28,29]), one may envisage many applications of these solutions in the future. In particular, one may apply the above solution for $N = 2$ to generate a new exactly solvable potential for the one-dimensional stationary Schrödinger equation. Other possibilities may appear if the technique used to derive the above-presented results is extended to involve the multivariable generalizations of the hypergeometric functions [30–33]. Finally, taking into account the relationship of the Fuchsian equations having five regular singularities with the Painlevé I-VI equations, applications to nonlinear problems are also possible [34,35].

**Author Contributions:** The authors contributed equally to this work. The authors read and approved the final version of the manuscript.

**Funding:** This research received no external funding.